\documentclass[11pt]{article}

\usepackage{amssymb}
\usepackage{latexsym}
\usepackage{latexsym}
\usepackage[tbtags]{amsmath}
\usepackage{epsfig}
\usepackage{amscd}
\usepackage{amsthm}
\usepackage{multirow}
\usepackage{float}
\usepackage{footmisc}
\usepackage{color}
\usepackage{enumerate}
\usepackage{bm}
\usepackage{wrapfig}
\usepackage{setspace}
\usepackage{adjustbox} % table size
\usepackage{xparse}
\usepackage{subfigure} % subfigures
\usepackage{tcolorbox}
\usepackage{lipsum}
\usepackage{blindtext}
\usepackage{mathtools}
\usepackage[titletoc,title]{appendix}
\usepackage[T1]{fontenc}
\usepackage{dsfont}
\usepackage{natbib}
\usepackage{algorithm}
\usepackage{algpseudocode}
\usepackage{url}
\usepackage[font=small,labelfont=bf]{caption}
\usepackage{authblk} % authors
\usepackage{hyperref}[] % hyperlink
\hypersetup{
	colorlinks = true,
	linkcolor = blue, %Colour of internal links
	filecolor = magenta,
	urlcolor = blue, %Colour for external hyperlinks
	citecolor = blue %Colour of citations
}

\usepackage{sectsty}

\newcommand{\mP}{\mathbb{P}}
\newcommand{\mE}{\mathbb{E}}

\newcommand{\mV}{\mathrm{Var}}

\newtheorem{lemma}{Lemma}

\allowdisplaybreaks

\begin{document}
	
	%\maketitle
	\begin{center}
		\LARGE \textbf{Comments on ``Testing Conditional Independence of Discrete Distributions''}
	\end{center}

	\vskip 1em
	\begin{center}
		Ilmun Kim \\
		Department of Statistics and Data Science \\
		%Department of Applied Statistics \\
		Yonsei University
	\end{center}
	\vskip 1em

	\begin{abstract}
		In this short note, we identify and address an error in the proof of Theorem 1.3 in \cite{canonne2018testing}, a recent breakthrough in conditional independence testing. After correcting the error, we show that the general sample complexity result established in \cite{canonne2018testing} remains the same. 
	\end{abstract}

	\section{A Closer Look at \cite{canonne2018testing}} \label{Section: Error}
	Let $(X,Y,Z)$ follow a distribution $p_{X,Y,Z}$ on a discrete domain $[\ell_1] \times [\ell_2] \times [n]$. Given $m$ i.i.d.~copies of $(X,Y,Z)$, the work of \cite{canonne2018testing} is concerned with the problem of testing conditional independence between $X$ and $Y$ given $Z$. Let $\mathcal{P}$ denote the set of all discrete distributions over $[\ell_1] \times [\ell_2] \times [n]$ and $\mathcal{P}_0 \subset \mathcal{P}$ denote the set of null distributions such that $X \perp\!\!\!\!\perp Y \,| \, Z$. Consider the set of alternative distributions, denoted by $\mathcal{P}_1(\varepsilon)$, which are $\varepsilon$ far way from $\mathcal{P}_0$ in the total variation distance:
	\begin{align*}
		\mathcal{P}_1(\varepsilon) = \Big\{ p \in \mathcal{P} : \inf_{q \in \mathcal{P}_0} \sup_{A \subseteq [\ell_1] \times [\ell_2] \times [n]} \big| p(A) - q(A) \big| \geq \varepsilon \Big\}.
	\end{align*} 
	Given these sets of distributions, the hypotheses of interest can be formulated as 
	\begin{align*}
		H_0: p_{X,Y,Z} \in \mathcal{P}_0 \quad \text{versus} \quad H_1: p_{X,Y,Z} \in \mathcal{P}_1(\varepsilon).
	\end{align*}
	For the above testing problem, Theorem 1.3 of \cite{canonne2018testing} proves that there exists a test with sample complexity
	\begin{align} \label{Eq: sample complexity}
		O\Biggl( \max \Biggl\{  \min \Biggl\{ \frac{n^{7/8}\ell_1^{1/4}\ell_2^{1/4} }{\varepsilon}, \ \frac{n^{6/7} \ell_1^{2/7}\ell_2^{2/7} }{\varepsilon^{8/7}}  \Biggr\}, \  \frac{n^{3/4}\ell_1^{1/2}\ell_2^{1/2} }{\varepsilon},  \ \frac{n^{2/3}\ell_1^{2/3}\ell_2^{1/3} }{\varepsilon^{4/3}}, \ \frac{n^{1/2}\ell_1^{1/2}\ell_2^{1/2} }{\varepsilon^{2}}  \Biggr\} \Biggr).
	\end{align}
	This result is achieved under Poissonization where the sample size $M$ follows a Poisson random variable with parameter $m$. The proof of Theorem 1.3 relies on several innovative lemmas including Claim 2.2:
	
	\vskip 1em 
	
	\noindent \textbf{Claim 2.2 of \cite{canonne2018testing}.} \emph{There exists an absolute constant $C>0$ such that, for $X \sim \mathrm{Poisson}(\lambda)$ and $a,b \geq 0$,}
	\begin{align*}
		\mV\bigl[ X \sqrt{\min(X,a) \min(X,b)} \mathds{1}(X \geq 4) \bigr] \leq C \mE\bigl[ X \sqrt{\min(X,a)\min(X,b)} \mathds{1}(X \geq 4) \bigr].
	\end{align*}
	Unfortunately, the upper bound in the claim turns out to be incorrect and therefore the proof of their Theorem 1.3 remains incomplete (confirmed via personal communication with one of the authors). We fix this error in Lemma~\ref{Lemma: Claim 2.2} below, which is the main technical contribution of this work. To illustrate the error in Claim 2.2, suppose that the parameter $\lambda$ is sufficiently large such that $X$ is approximately $N(\lambda,\lambda)$ by the central limit theorem. Moreover, by taking $a$ and $b$ relatively small in comparison to $\lambda$, one can have an approximation
	\begin{align*}
		X \sqrt{\min(X,a) \min(X,b)} \mathds{1}(X \geq 4) \approx N(\lambda,\lambda) \sqrt{ab}. 
	\end{align*}
	Hence the variance and the expectation in Claim 2.2 approximately become $\lambda ab$ and $\lambda \sqrt{ab}$, respectively, which contradicts the given inequality when $a,b \geq 1$. 
	
	\section{Main Result} 
	\sloppy We now show that the upper bound in Claim 2.2 holds if we add an additional factor of $\max\{ab,\sqrt{a}b, \sqrt{ab}\}$ for $b \geq a \geq 0$ in the upper bound. 
		
	\begin{lemma}[Correction of Claim 2.2] \label{Lemma: Claim 2.2}
		There exists an absolute constant $C>0$ such that, for $X \sim \mathrm{Poisson}(\lambda)$ and $b \geq a \geq 0$,
		\begin{align*}
		& \mV\bigl[ X \sqrt{\min(X,a) \min(X,b)} \mathds{1}(X \geq 4) \bigr]\\[.5em] 
		\leq~ &  C  \max\{ab,\sqrt{a}b, \sqrt{ab}\} \mE\bigl[ X \sqrt{\min(X,a)\min(X,b)} \mathds{1}(X \geq 4) \bigr].
		\end{align*}
	\end{lemma}	
	\noindent %It is worth pointing out that the above lemma further assumes $a,b \geq 2$. This restriction is not an issue in its application where $a$ and $b$ correspond to the domain sizes $\ell_1 \geq 2$ and $\ell_2 \geq 2$, respectively. We also note that one can improve the factor of $ab$ in the upper bound with more effort. Nevertheless, the current form is sharp enough to reprove the sample complexity~(\ref{Eq: sample complexity}). Due to the non-linearity of $X \sqrt{\min(X,a) \min(X,b)} \mathds{1}(X \geq 4)$ in $X$, the proof of Lemma~\ref{Lemma: Claim 2.2} turns out to be non-trivial, requiring a delicate case-by-case analysis. The details can be found in Section~\ref{Section: Proof}. 
	It is worth pointing out that one can improve the factor of $\max\{ab,\sqrt{a}b,\sqrt{ab}\}$ in the upper bound with more effort. Nevertheless, the current form is sharp enough to reprove the sample complexity~(\ref{Eq: sample complexity}). We also note that in its application, $a$ and $b$ correspond to the domain sizes $\ell_1 \geq 2$ and $\ell_2 \geq 2$, respectively, in which case the additional factor becomes $ \max\{ab,\sqrt{a}b,\sqrt{ab}\}  = ab$. Due to the non-linearity of $X \sqrt{\min(X,a) \min(X,b)} \mathds{1}(X \geq 4)$ in $X$, the proof of Lemma~\ref{Lemma: Claim 2.2} turns out to be non-trivial, requiring a delicate case-by-case analysis. The details can be found in Section~\ref{Section: Proof}. 
	
	\vskip 1em
	
	Having provided a correction of Claim 2.2, we next reprove Theorem 1.3 of \cite{canonne2018testing}. There is essentially one place where Claim 2.2 is used in \cite{canonne2018testing}, more specifically in their Lemma 5.3. Once we establish Lemma 5.3, the rest of the proof remains the same.  Before stating the result, we need to introduce some notation. For simplicity, we denote the conditional distribution of $X$ and $Y$ given $Z=z$ by $p_z$ and the product distribution of their marginals by $q_z$. Define $\varepsilon_z = \sup_{A \subseteq [n]} |p_z(A) - q_z(A)|$ and write $\varepsilon_z' = \varepsilon_z / \sqrt{4\ell_1\ell_2}$. We further let $\omega_z = \sqrt{\min(\sigma_z,\ell_1)\min(\sigma_z,\ell_2)}$. With this notation in place, we reprove Lemma 5.3 of \cite{canonne2018testing} as follows.
	
	\begin{lemma}[Lemma 5.3 of \cite{canonne2018testing}] \label{Lemma: Lemma 5.3}
		Consider $D = \sum_{z \in [n]} \sigma_z \omega_z \varepsilon_z^{\prime 2} \mathds{1}(\sigma_z \geq 4)$ where $\sigma_1,\ldots,\sigma_n$ are independent Poisson random variables with the corresponding parameters $m \mP(Z=1), \ldots, m \mP(Z=n)$. Then there exists a constant $C>0$ such that $\mV[D] \leq C \mE[D]$.
	\begin{proof}
		Using the independence of $\sigma_1,\ldots,\sigma_n$, we have the identity
		\begin{align*}
			\mV[D] ~=~ \sum_{z \in [n]}  \varepsilon_z^{\prime 4} \mV[ \sigma_z \omega_z \mathds{1}(\sigma_z \geq 4)].
		\end{align*}
		Moreover, Lemma~\ref{Lemma: Claim 2.2} along with a trivial bound $\varepsilon_z' \leq 1/\sqrt{4 \ell_1\ell_2}$ yields
		\begin{align*}
			 \sum_{z \in [n]}  \varepsilon_z^{\prime 4} \mV[ \sigma_z \omega_z \mathds{1}(\sigma_z \geq 4)] ~\leq~& C\sum_{z \in [n]}  \varepsilon_z^{\prime 4} \ell_1 \ell_2 \mE[ \sigma_z \omega_z \mathds{1}(\sigma_z \geq 4)] \\[.5em]
			 \leq ~ & C \sum_{z \in [n]}  \varepsilon_z^{\prime 2} \mE[ \sigma_z \omega_z \mathds{1}(\sigma_z \geq 4)] = C \mE[D].
		\end{align*}
		This completes the proof of Lemma~\ref{Lemma: Lemma 5.3}, which in turn proves Theorem 1.3 of \cite{canonne2018testing}.  
	\end{proof}
	\end{lemma}

	 \section{Proof of Lemma~\ref{Lemma: Claim 2.2}}  \label{Section: Proof}
	 It now remains to prove Lemma~\ref{Lemma: Claim 2.2}. We begin with the notation used in the proof. 
	 
	 \vskip .5em
	 
	 \noindent \textbf{Notation.} We use $C,C'$ to denote positive absolute constants whose values may vary from line to line. For constants $x$ and $y$, we write $x \lesssim y$ (resp.~$x \gtrsim y$) if there exists an absolute positive constant $C$ (resp.~$C'$) such that $x \leq C y$ (resp.~$x \geq C'y$). $a \wedge b$ indicates $\min(a,b)$. For simplicity, we often write 
	 \begin{align*}
	 	X_{a,b} := X \sqrt{(X \wedge a) (X \wedge b)} \mathds{1}(X \geq 4).
	 \end{align*}
	 
	 \vskip .5em
	 
	 \noindent Our proof relies on two other claims in \cite{canonne2018testing}. We recall them here for completeness.
	 
	 \vskip .5em
	 
	 \noindent \textbf{Claim 2.1 of \cite{canonne2018testing}.} \emph{For $X \sim \mathrm{Poisson}(\lambda)$, the following inequality holds}
	 \begin{align} \label{Eq: Claim 2.1}
	 	\mV[X \mathds{1}(X \geq 4)] ~\lesssim~ \mE[X \mathds{1}(X \geq 4)].
	 \end{align}
	 
	 \noindent \textbf{Claim 2.3 of \cite{canonne2018testing}} \emph{For $X \sim \mathrm{Poisson}(\lambda)$ and integers $a,b \geq 2$, the following inequality holds}
	 \begin{align} \label{Eq: Claim 2.3}
	 	\mE\bigl[X \sqrt{(X \wedge a) (X \wedge b)} \mathds{1}(X \geq 4)\bigr] ~\gtrsim~ \bigl( \lambda\sqrt{(\lambda \wedge a) (\lambda \wedge b)}  \wedge \lambda^4 \bigr).
	 \end{align}
	 
	 \vskip .5em
	 
	 \noindent We split the proof into two parts depending on whether $a \geq 2$ or $0 \leq a < 2$.
	 
	 \subsection{Case where $a \geq 2$} \label{Section: Case 1}
	 We start with the first case where $a \geq 2$. In our analysis, it is more convenient to work with an alternative form of the variance of $X_{a,b}$. First notice that for any random variable $W$ and its i.i.d.~copy $W'$, the variance of $W$ is equal to
	 \begin{align*}
	 	\mV[W] = \frac{\mE\bigl[ (W - W')^2 \bigr]}{2}.
	 \end{align*}
	 Let $Y$ be an i.i.d.~copy of $X \sim \mathrm{Poisson}(\lambda)$. Using the alternative formula of the variance, we have the identity given as
	 \begin{align*}
	 	 2\mathrm{Var}[X_{a,b}] 
	 	=  \mE\bigl[ \bigl(X \sqrt{(X \wedge a) (X \wedge b)} \mathds{1}(X \geq 4) - Y \sqrt{(Y \wedge a) (Y \wedge b)} \mathds{1}(Y \geq 4) \bigr)^2 \bigr].
	 \end{align*}
	 Equivalently, by letting $p$ be the probability mass function (pmf) of $\mathrm{Poisson}(\lambda)$, 
	 \begin{align*}
	 	2\mathrm{Var}[X_{a,b}] ~=~& \underbrace{\sum_{x=4}^\infty \sum_{y=4}^\infty \Bigl( x \sqrt{(x \wedge a)(x \wedge b)} - y \sqrt{(y \wedge a)(y \wedge b)} \Bigr)^2 p(x) p(y)}_{:=(\mathrm{I})} \\[.5em]
	 	+ ~ &\underbrace{\sum_{x=4}^\infty \sum_{y=0}^3 \Bigl( x \sqrt{(x \wedge a)(x \wedge b)} \Bigr)^2 p(x) p(y)}_{:=(\mathrm{II})} \\[.5em]
	 	+ ~ & \underbrace{\sum_{x=0}^3 \sum_{y=4}^\infty \Bigl( y \sqrt{(y \wedge a)(y \wedge b)} \Bigr)^2 p(x) p(y)}_{:=(\mathrm{III})},
	 \end{align*}
	 where we use the fact that both $X_{a,b}$ and $Y_{a,b}$ are equal to zero when $x \leq 3$ and $y \leq 3$. 
	 
	 \vskip 1em
	 
	 \noindent \textbf{Analysis of (II) and (III).} Let us start bounding the second term~(II) as it is simpler to analyze than the first term~(I). The third term can be handled similarly by symmetry. We proceed by considering the two cases (i)~$\lambda \geq 1 $ and (ii)~$\lambda < 1$, separately.
	 
	 \begin{itemize}
	 	\item Suppose that $\lambda \geq 1$. Using the explicit form of the Poisson pmf and its fourth moment expression, we have 
	 	\begin{align*}
	 		\sum_{x=4}^\infty \sum_{y=0}^3 \Bigl( x \sqrt{(x \wedge a)(x \wedge b)} \Bigr)^2 p(x)p(y)  ~\leq~& \sum_{x=4}^\infty x^4 p(x)  \sum_{y=0}^3 p(y) 
	 		\\[.5em]
	 		\leq ~  & \mE[X^4] \biggl(1 + \lambda + \frac{\lambda^2}{2} + \frac{\lambda^3}{6}\biggr) e^{-\lambda} \\[.5em]
	 		~=~ & (\lambda^4 + 6 \lambda^3 + 7 \lambda^2 + \lambda)  \biggl(1 + \lambda + \frac{\lambda^2}{2} + \frac{\lambda^3}{6}\biggr) e^{-\lambda}.
	 	\end{align*}
	 	We would like to verify that the above bound is less than or equal to $\mE[X_{a,b}]$, up to a constant. To this end, recall Claim 2.3 of \cite{canonne2018testing} in (\ref{Eq: Claim 2.3}):
	 	\begin{align*}
	 		\mE[X_{a,b}] ~\gtrsim~ \bigl( \lambda\sqrt{(\lambda \wedge a) (\lambda \wedge b)}  \wedge \lambda^4 \bigr).
	 	\end{align*}
 		Since we assume $a,b \geq 2$, we have $\sqrt{(\lambda \wedge a) (\lambda \wedge b)} \geq 1$ and so
	 	\begin{align*}
	 		\mE[X_{a,b}] ~\gtrsim~ &  \bigl( \lambda \wedge \lambda^4 \bigr) \\[.5em]
	 		\gtrsim ~ & (\lambda^4 + 6 \lambda^3 + 7 \lambda^2 + \lambda) \biggl(1 + \lambda + \frac{\lambda^2}{2} + \frac{\lambda^3}{6}\biggr)  e^{-\lambda},
	 	\end{align*}
	 	where the last inequality uses the fact that the following function
	 	\begin{align*}
	 		h(\lambda) := \frac{\lambda\wedge \lambda^4}{(\lambda^4 + 6 \lambda^3 + 7 \lambda^2 + \lambda)  \Bigl(1 + \lambda + \frac{\lambda^2}{2} + \frac{\lambda^3}{6} \Bigr) e^{-\lambda}}
	 	\end{align*}
	 	is bounded below by some positive constant (numerically $\approx 0.0119$) for all $\lambda \geq 1$. Therefore when $\lambda \geq 1$,
	 	\begin{align*}
	 		\sum_{x=4}^\infty \sum_{y=0}^3 \Bigl( x \sqrt{(x \wedge a)(x \wedge b)} \Bigr)^2 p(x) p(y) ~\lesssim~ \mE[X_{a,b}].
	 	\end{align*}
 		\item Now consider the case where $\lambda <1$. In this case, we observe $\mE[X_{a,b}] \gtrsim \lambda^4$. Also note that
 		\begin{align*}
 			& \sum_{x=4}^\infty \sum_{y=0}^3 \Bigl( x \sqrt{(x \wedge a)(x \wedge b)} \Bigr)^2 p(x) p(y) \\[.5em] ~\overset{\mathrm{(i)}}{\leq}~ & \sum_{x=4}^\infty \Bigl( x \sqrt{(x \wedge a)(x \wedge b)} \Bigr)^2 p(x) \\[.5em]
 			\overset{\mathrm{(ii)}}{\leq} ~ & \sqrt{ab} \sum_{x=4}^\infty x^2  \sqrt{(x \wedge a)(x \wedge b)} \frac{\lambda^x e^{-\lambda}}{x!} \\[.5em]
 			= ~ &   \sqrt{ab} \lambda \sum_{x=4}^\infty x  \sqrt{(x \wedge a)(x \wedge b)} \frac{\lambda^{x-1} e^{-\lambda}}{(x-1)!} \\[.5em]
 			\overset{\mathrm{(iii)}}{=} ~ &  \sqrt{ab} \lambda \sum_{x=3}^\infty (x+1)  \sqrt{((x+1) \wedge a)((x+1) \wedge b)} \frac{\lambda^{x} e^{-\lambda}}{x!} \\[.5em]
 			\overset{\mathrm{(iv)}}{\leq} ~ &  4 \sqrt{ab} \lambda \sum_{x=3}^\infty x \sqrt{(x \wedge a)(x \wedge b)} \frac{\lambda^{x} e^{-\lambda}}{x!} \\[.5em]
 			= ~ & 4 \sqrt{ab} \lambda \times 3 \sqrt{(3 \wedge a) (3 \wedge b)} \frac{\lambda^3 e^{-\lambda}}{3!} + 4 \sqrt{ab} \lambda \sum_{x=4}^\infty x \sqrt{(x \wedge a)(x \wedge b)} \frac{\lambda^{x} e^{-\lambda}}{x!} \\[.5em]
 			\lesssim ~ & \sqrt{ab} \lambda^4 + \sqrt{ab} \mE[X_{a,b}],
 		\end{align*}
 		where step~(i) uses the fact $p$ is a probability mass function, step~(ii) uses the basic inequality $\sqrt{(x \wedge a) (y \wedge b)} \leq \sqrt{ab}$, step~(iii) follows by change of variables and step~(iv) follows since $x +1 \leq 2x$ for all $x \geq 3$. Therefore we can conclude that 
 		\begin{align*}
 			\sum_{x=4}^\infty \sum_{y=0}^3 \Bigl( x \sqrt{(x \wedge a)(x \wedge b)} \Bigr)^2 p(x) p(y) ~\lesssim~ \sqrt{ab} \mE[X_{a,b}].
 		\end{align*}
 		Similarly, 
 		\begin{align*}
 			\sum_{x=0}^3 \sum_{y=4}^\infty \Bigl( y \sqrt{(y \wedge a)(y \wedge b)} \Bigr)^2 p(x) p(y) ~\lesssim~ \sqrt{ab} \mE[X_{a,b}].
 		\end{align*}
	 \end{itemize}
	
	 \vskip 1em
	 
	 \noindent \textbf{Analysis of (I).} Now, we handle the first term (I) in the variance expansion:
	 \begin{align*}
	 	(\mathrm{I}) ~=~ \sum_{x=4}^\infty \sum_{y=4}^\infty \underbrace{\Bigl( x \sqrt{(x \wedge a)(x \wedge b)} - y \sqrt{(y \wedge a)(y \wedge b)} \Bigr)^2 p(x) p(y)}_{:=u_{x,y,a,b}}.
	 \end{align*}
 	We separate the cases into three: (i)~$b\geq a \geq 4$, (ii)~$2 \leq a,b \leq 4$ and (iii)~$2 \leq a < 4 \leq b$ and proceed our analysis.
 	\begin{itemize}
 		\item We first assume that (i)~$b \geq a \geq 4$ and decompose the above double summations into 9 terms:
 		\begin{align*}
 			\sum_{x=4}^\infty \sum_{y=4}^\infty u_{x,y,a,b} ~ = ~ & \sum_{x=4}^{a} \sum_{y=4}^{a} u_{x,y,a,b} +  \sum_{x=4}^{a} \sum_{y=a+1}^{b} u_{x,y,a,b} +  \sum_{x=4}^{a} \sum_{y=b+1}^{\infty} u_{x,y,a,b} \\[.5em]  
 			+ ~ & \sum_{x=a+1}^{b} \sum_{y=4}^{a} u_{x,y,a,b} +  \sum_{x=a+1}^{b}\sum_{y=a+1}^{b} u_{x,y,a,b}  +   \sum_{x=a+1}^{b} \sum_{y=b+1}^{\infty} u_{x,y,a,b} \\[.5em]  
 			+ ~ & \sum_{x=b+1}^{\infty} \sum_{y=4}^{a} u_{x,y,a,b} +  \sum_{x=b+1}^{\infty} \sum_{y=a+1}^{b} u_{x,y,a,b} +  \sum_{x=b+1}^{\infty}  \sum_{y=b+1}^{\infty} u_{x,y,a,b} \\[.5em]
 			= ~ & \mathrm{(I)}_1 + \cdots + \mathrm{(I)}_9.
 		\end{align*}
 		Starting with the first term~$\mathrm{(I)}_1$,
 		\begin{align*}
 			\sum_{x=4}^{a} \sum_{y=4}^{a} u_{x,y,a,b} ~=~ & \sum_{x=4}^{a} \sum_{y=4}^{a} \Bigl( x \sqrt{(x \wedge a)(x \wedge b)} - y \sqrt{(y \wedge a)(y \wedge b)} \Bigr)^2 p(x) p(y) \\[.5em]
 			= ~ & \sum_{x=4}^{a} \sum_{y=4}^{a} (x^2 - y^2)^2 p(x)p(y) \\[.5em]
 			= ~ &  \sum_{x=4}^{a} \sum_{y=4}^{a} (x - y)^2 (x+y)^2 p(x)p(y) \\[.5em] 
 			\leq ~ & 4a^2 \sum_{x=4}^{a} \sum_{y=4}^{a} (x - y)^2 p(x)p(y) \\[.5em]
 			\lesssim ~ & a^2 \mV[X \mathds{1}(X \geq 4)],
 		\end{align*}
 		where the last inequality follows by the fact that for i.i.d.~random variables $X$ and $Y$,
 		\begin{align*}
 			& \mV[X \mathds{1}(X \geq 4)]  \\[.5em]
 			~=~& \frac{1}{2} \mE\bigl[ \bigl(X \mathds{1}(X \geq 4) - Y \mathds{1}(Y \geq 4)  \bigr)^2  \bigr] \\[.5em]
 			= ~ &  \frac{1}{2}\sum_{x=4}^{a} \sum_{y=4}^{a} (x - y)^2 p(x)p(y) + \frac{1}{2} \sum_{x=4}^{a} \sum_{y=0}^{3} x^2 p(x)p(y) + \frac{1}{2} \sum_{x=0}^{3} \sum_{y=4}^{a} y^2 p(x)p(y).
 		\end{align*}
 		For the second term~$\mathrm{(I)}_2$,
 		\begin{align*}
 			\sum_{x=4}^{a} \sum_{y=a+1}^{b} u_{x,y,a,b} ~=~& \sum_{x=4}^{a} \sum_{y=a+1}^{b}  \Bigl( x \sqrt{(x \wedge a)(x \wedge b)} - y \sqrt{(y \wedge a)(y \wedge b)} \Bigr)^2 p(x) p(y) \\[.5em]
 			=~& \sum_{x=4}^{a} \sum_{y=a+1}^{b}  \bigl( x^2 - y \sqrt{a y} \bigr)^2 p(x) p(y). 
 		\end{align*}
 		Since $x \leq a < y \leq b$, we have $y\sqrt{y} \sqrt{a} - x^2 \geq 0$ and 
 		\begin{align*}
 			y\sqrt{y} \sqrt{a} -x^2 ~=~ & x \sqrt{y} \sqrt{a} - x^2  - x \sqrt{y} \sqrt{a} + y \sqrt{y} \sqrt{a} \\[.5em]
 			= ~ & x(\sqrt{y}\sqrt{a} - x) + \sqrt{y} \sqrt{a} (y-x) \\[.5em]
 			\leq ~ & x(y- x) +  \sqrt{y}\sqrt{a} (y-x) \\[.5em] 
 			=~ & (x + \sqrt{y}\sqrt{a})(y-x), 
 		\end{align*}
 		which implies 
 		\begin{align*}
 			( x^2 - y\sqrt{y} \sqrt{a} )^2 ~\leq~  (x + \sqrt{y}\sqrt{a})^2(y-x)^2 \leq 4ab (y-x)^2.
 		\end{align*}
 		Using this inequality, we have
 		\begin{align*}
 			\sum_{x=4}^{a} \sum_{y=a+1}^{b} u_{x,y,a,b}  ~\lesssim~ ab  \mV[X \mathds{1}(X \geq 4)].
 		\end{align*}
 		For the third term~$\mathrm{(I)}_3$, observe that
 		\begin{align*}
 			\sum_{x=4}^{a} \sum_{y=b+1}^{\infty} u_{x,y,a,b} ~=~ &  \sum_{x=4}^{a} \sum_{y=b+1}^{\infty}   \Bigl( x \sqrt{(x \wedge a)(x \wedge b)} - y \sqrt{(y \wedge a)(y \wedge b)} \Bigr)^2 p(x) p(y) \\[.5em]
 			= ~ &  \sum_{x=4}^{a} \sum_{y=b+1}^{\infty}  \bigl( x^2 - y \sqrt{ab} \bigr)^2 p(x) p(y).
 		\end{align*}
 		In this case, we have $x \leq a \leq b < y$ and
 		\begin{align*}
 			0 ~\leq~&  y\sqrt{ab} - x^2 \\[.5em]
 			= ~ & x\sqrt{ab} - x^2 + y \sqrt{ab} - x\sqrt{ab}  \\[.5em]
 			= ~ & x(\sqrt{ab} - x) + \sqrt{ab}(y - x) \\[.5em]
 			\leq ~ & (x+\sqrt{ab}) (y-x),
 		\end{align*}
 		where the last inequality holds since $\sqrt{ab} < y$, which in turn implies that
 		\begin{align*}
 			\sum_{x=4}^{a} \sum_{y=b+1}^{\infty} u_{x,y,a,b} ~\lesssim~ ab  \mV[X \mathds{1}(X \geq 4)].
 		\end{align*}
 		By symmetry, the fourth term~$\mathrm{(I)}_4$ can be similarly analyzed as the second term. For the fifth term~$\mathrm{(I)}_5$,
 		\begin{align*}
 			\sum_{x=a+1}^{b}\sum_{y=a+1}^{b} u_{x,y,a,b}  ~=~& \sum_{x=a+1}^{b}\sum_{y=a+1}^{b} \Bigl( x \sqrt{(x \wedge a)(x \wedge b)} - y \sqrt{(y \wedge a)(y \wedge b)} \Bigr)^2 p(x) p(y) \\[.5em]
 			= ~ &  \sum_{x=a+1}^{b}\sum_{y=a+1}^{b} \bigl( x \sqrt{a x} - y \sqrt{ay} \bigr)^2 p(x) p(y) \\[.5em]
 			= ~ &  a \sum_{x=a+1}^{b}\sum_{y=a+1}^{b} \bigl( x \sqrt{x} - y \sqrt{y} \bigr)^2 p(x) p(y).
 		\end{align*} 
 		Without loss of generality, assume $a < x \leq y \leq b$ and note that
 		\begin{align*}
 			0 ~ \leq ~ & y \sqrt{y} - x\sqrt{x}\\[.5em]
 			= ~ & x \sqrt{b} - x\sqrt{x} + y\sqrt{b} - x \sqrt{b} + y \sqrt{y} - y \sqrt{b}  \\[.5em]
 			= ~ & x (\sqrt{b} - \sqrt{x}) + \sqrt{b}(y-x) + y (\sqrt{y} - \sqrt{b}) \\[.5em]
 			\leq ~ & y (\sqrt{b} - \sqrt{x}) + \sqrt{b}(y-x) + y (\sqrt{y} - \sqrt{b}) \\[.5em]
 			= ~ & y (\sqrt{y} - \sqrt{x}) + \sqrt{b} (y - x) \\[.5em]
 			= ~ & y \frac{y-x}{\sqrt{y} + \sqrt{x}} + \sqrt{b} (y-x).
 		\end{align*}
 		Therefore we have
 		\begin{equation}
			\begin{aligned}  \label{Eq: analysis}
				(x\sqrt{x} - y \sqrt{y})^2 ~\leq~&  \frac{2y^2}{(\sqrt{x} + \sqrt{y})^2} (x-y)^2 + 2b(x-y)^2\\[.5em] 
				\leq ~ & 2 y (x-y)^2 + 2b(x-y)^2,
			\end{aligned}
 		\end{equation}
 		which implies along with condition $y \leq b$ that
 		\begin{align*}
 			\sum_{x=a+1}^{b}\sum_{y=a+1}^{b} u_{x,y,a,b}  ~\lesssim~  ab  \mV[X \mathds{1}(X \geq 4)].
 		\end{align*}
 		For the sixth term~$\mathrm{(I)}_6$, it can be seen that 
 		\begin{align*}
 			\sum_{x=a+1}^{b} \sum_{y=b+1}^{\infty} u_{x,y,a,b} ~=~ & \sum_{x=a+1}^{b} \sum_{y=b+1}^{\infty}\Bigl( x \sqrt{(x \wedge a)(x \wedge b)} - y \sqrt{(y \wedge a)(y \wedge b)} \Bigr)^2 p(x) p(y) \\[.5em]
 			= ~ & \sum_{x=a+1}^{b} \sum_{y=b+1}^{\infty}\bigl( x \sqrt{a x} - y \sqrt{ab} \bigr)^2 p(x) p(y) \\[.5em]
 			= ~ & a \sum_{x=a+1}^{b} \sum_{y=b+1}^{\infty}\bigl( x \sqrt{x} - y \sqrt{b} \bigr)^2 p(x) p(y).
 		\end{align*}
 		Since $a < x \leq b < y$, we have
 		\begin{align*}
 			0 ~<~ & y \sqrt{b}  - x \sqrt{x}~ = \sqrt{x}(\sqrt{x}\sqrt{b} - x) + \sqrt{b}(y-x) \\[.5em]
 			\leq~ & \sqrt{b}(\sqrt{x}\sqrt{b} - x) + \sqrt{b}(y-x) ~\leq~ 2\sqrt{b}(y -x),
 		\end{align*}
 		which implies 
 		\begin{align*}
 			\sum_{x=a+1}^{b} \sum_{y=b+1}^{\infty} u_{x,y,a,b} ~\lesssim~ ab  \mV[X \mathds{1}(X \geq 4)].
 		\end{align*}
 		The seventh term~$\mathrm{(I)}_7$ and the eighth term~$\mathrm{(I)}_8$ can be analyzed as the third term~$\mathrm{(I)}_3$ and the sixth term~$\mathrm{(I)}_6$, respectively. The final term~$\mathrm{(I)}_9$ is 
 		\begin{align*}
 			\sum_{x=b+1}^{\infty}  \sum_{y=b+1}^{\infty} u_{x,y,a,b} ~=~ &\sum_{x=b+1}^{\infty}  \sum_{y=b+1}^{\infty} \Bigl( x \sqrt{(x \wedge a)(x \wedge b)} - y \sqrt{(y \wedge a)(y \wedge b)} \Bigr)^2 p(x) p(y) \\[.5em]
 			= ~ & \sum_{x=b+1}^{\infty}  \sum_{y=b+1}^{\infty} \bigl( x \sqrt{ab} - y \sqrt{ab} \bigr)^2 p(x) p(y) \\[.5em]
 			= ~ & ab \sum_{x=b+1}^{\infty}  \sum_{y=b+1}^{\infty} ( x - y )^2 p(x) p(y)\\[.5em]
 			\lesssim ~ & ab  \mV[X \mathds{1}(X \geq 4)].
 		\end{align*}
 		In summary, we have established that
 		\begin{align*}
 			(\mathrm{I}) ~=~ \sum_{x=4}^\infty \sum_{y=4}^\infty u_{x,y,a,b} ~\lesssim~ & ab  \mV[X \mathds{1}(X \geq 4)] ~	\overset{\mathrm{(i)}}{\lesssim}  ~ ab \mE[X \mathds{1}(X \geq 4)] \\[.5em]
 			\lesssim ~ & ab \mE[X \sqrt{(X \wedge a)} \sqrt{(X \wedge b)} \mathds{1}(X \geq 4)], 
 		\end{align*}
 		where step~(i) uses Claim 2.1 of \cite{canonne2018testing} in (\ref{Eq: Claim 2.1}). Combining the pieces yields
 		\begin{align*}
 			\mV[X_{a,b}] ~\lesssim~ ab \mE[X_{a,b}],
 		\end{align*}
 		as desired. 
 		\item We are only left with the cases where (ii)~$2 \leq a,b < 4$ and (iii)~$2 \leq a < 4 \leq b$, and only need to re-analyze the term~(I). Since both $x,y \geq 4$, it holds that $ \sqrt{(x \wedge a)(x \wedge b)} = \sqrt{ab}$ under case~(i) $2 \leq a,b <4$, and therefore we have
 		\begin{align*}
 			\sum_{x=4}^\infty \sum_{y=4}^\infty u_{x,y,a,b} ~=~ & \sum_{x=4}^\infty \sum_{y=4}^\infty \Bigl( x \sqrt{(x \wedge a)(x \wedge b)} - y \sqrt{(y \wedge a)(y \wedge b)} \Bigr)^2 p(x) p(y) \\[.5em]
 			= ~ & ab  \sum_{x=4}^\infty \sum_{y=4}^\infty ( x - y)^2 p(x) p(y) \\[.5em]
 			\lesssim ~ & ab  \mV[X \mathds{1}(X \geq 4)] ~ \lesssim ~ ab \mE[X_{a,b}]. 
 		\end{align*}
 		Similarly, when case (ii) $2 \leq a < 4 \leq b$ holds, 
 		\begin{align*}
 			\sum_{x=4}^\infty \sum_{y=4}^\infty u_{x,y,a,b} ~=~ &  a  \sum_{x=4}^\infty \sum_{y=4}^\infty \underbrace{\Bigl( x \sqrt{x \wedge b} - y \sqrt{y \wedge b} \Bigr)^2 p(x)p(y)}_{:=v_{x,y,b}} \\[.5em]
 			= ~ & a \sum_{x=4}^b \sum_{y=4}^b v_{x,y,b} + a \sum_{x=b+1}^\infty \sum_{y=4}^b v_{x,y,b} \\[.5em] 
 			+ ~ & a \sum_{x=4}^b \sum_{y=b+1}^\infty v_{x,y,b} + a \sum_{x=b+1}^\infty \sum_{y=b+1}^\infty v_{x,y,b} \\[.5em]
 			= ~ & (\mathrm{I})_{1'} + \cdots + (\mathrm{I})_{4'}. 
 		\end{align*}
 		Following a similar analysis in (\ref{Eq: analysis}), the first term $(\mathrm{I})_{1'}$ in the decomposition satisfies
 		\begin{align*}
 			a \sum_{x=4}^b \sum_{y=4}^b v_{x,y,b} ~\lesssim~ ab \sum_{x=4}^b \sum_{y=4}^b (x-y)^2 p(x) p(y) ~\lesssim~ ab \mE[X_{a,b}].
 		\end{align*}
 		The second term $(\mathrm{I})_{2'}$ (and also the third term $(\mathrm{I})_{3'}$ by symmetry) satisfies 
 		\begin{align*}
 			a \sum_{x=b+1}^\infty \sum_{y=4}^b v_{x,y,b} ~=~ & a \sum_{x=b+1}^\infty \sum_{y=4}^b (x \sqrt{b} - y \sqrt{y})^2 p(x)p(y) \\[.5em]
 			\overset{(\star)}{\lesssim} ~ & ab \sum_{x=b+1}^\infty \sum_{y=4}^b (x - y)^2 p(x)p(y) ~\lesssim~ ab \mE[X_{a,b}],
 		\end{align*}
 		where inequality~($\star$) follows since it holds 
 		\begin{align*}
 			0 \leq x \sqrt{b} - y\sqrt{y} = (x-y) \sqrt{b} + \sqrt{y}(\sqrt{yb} - y) \leq  2 \sqrt{b}(x-y) \quad \text{when $4 \leq y \leq b < x$.}
 		\end{align*}
 		For the fourth term~$(\mathrm{I})_{4'}$, we have
 		\begin{align*}
 			a \sum_{x=b+1}^\infty \sum_{y=b+1}^\infty v_{x,y,b} ~=~  & ab \sum_{x=b+1}^\infty \sum_{y=b+1}^\infty (x-y)^2 p(x)p(y)  \\[.5em]
 			\lesssim~  & ab \mE[X_{a,b}].
 		\end{align*}
 		Combining the pieces yields
 		\begin{align*}
 			\sum_{x=4}^\infty \sum_{y=4}^\infty u_{x,y,a,b}  ~\lesssim~  ab \mE[X_{a,b}].
 		\end{align*}
 	\end{itemize}
	 This completes the proof of Lemma~\ref{Lemma: Claim 2.2} when $a \geq 2$.

	 \subsection{Case where $0 \leq a < 2$}
	 Next we consider the remaining case where $0 \leq a < 2$. In this case along with $b \geq 2$, observe that $X_{a,b} = X \sqrt{a (X \wedge b)} \mathds{1}(X \geq 4)$. Consequently it follows that
	 \begin{align*}
	 	\mV[X_{a,b}] ~=~ & a \mV\bigl[X \sqrt{(X \wedge b)} \mathds{1}(X \geq 4)\bigr] \\[.5em]
	 	= ~ & \frac{a}{2}\mV\bigl[X \sqrt{(X \wedge 2)(X \wedge b)} \mathds{1}(X \geq 4)\bigr] \\[.5em]
	 	\lesssim ~ & ab \mE\bigl[ X \sqrt{(X \wedge 2)(X \wedge b)} \mathds{1}(X \geq 4) \bigr],
	 \end{align*}
	 where the last inequality uses the previously established result in Section~\ref{Section: Case 1} for $a=2$ and $b \geq 2$. The desired result then follows by noting that 
	 \begin{align*}
	 	ab \mE\bigl[ X \sqrt{(X \wedge 2)(X \wedge b)} \mathds{1}(X \geq 4) \bigr] ~=~ \sqrt{a}b \mE\bigl[ X \sqrt{(X \wedge a)(X \wedge b)} \mathds{1}(X \geq 4) \bigr].
	 \end{align*}
 	When $0 \leq a,b <2$, Claim 2.1 of \cite{canonne2018testing} in (\ref{Eq: Claim 2.1}) yields
 	\begin{align*}
 		\mV[X_{a,b}] ~=~ & ab \mV\bigl[X \mathds{1}(X \geq 4)\bigr] \\[.5em]
 		\lesssim ~ & ab \mE\bigl[X \mathds{1}(X \geq 4)\bigr] \\[.5em] 
 		\lesssim ~ & \sqrt{ab} \mE\bigl[ X \sqrt{(X \wedge a)(X \wedge b)} \mathds{1}(X \geq 4) \bigr].
 	\end{align*}
	 Therefore the statement of Lemma~\ref{Lemma: Claim 2.2} holds when $0 \leq a < 2$. Combining this with the previous result in Section~\ref{Section: Case 1} completes the proof of Lemma~\ref{Lemma: Claim 2.2}.
	
	\paragraph{Acknowledgements.} The author would like to thank Cl{\'e}ment L.~Canonne for his careful proofreading and helpful comments. The author also thanks Matey Neykov, Sivaraman Balakrishnan and Larry Wasserman for their helpful discussion. Finally, the author acknowledges support from the Yonsei University Research Fund of 2021-22-0332 as well as support from the Basic Science Research Program through the National Research Foundation of Korea (NRF) funded by the Ministry of Education (2022R1A4A1033384).
	\bibliographystyle{apalike}
	\bibliography{reference}
\end{document}